\documentclass [12pt]{amsart}

\usepackage[top=100pt,bottom=60pt,left=96pt,right=92pt]{geometry}
\usepackage[utf8]{inputenc}
\usepackage[T1]{fontenc}
\usepackage[catalan,german,spanish,dutch,british]{babel}
\usepackage{lmodern}
\usepackage[pdftex]{graphicx}
\usepackage{tikz}
\usetikzlibrary{matrix}
\usepackage{pinlabel}
\usepackage{faktor} 

\usepackage{mathtools} 

\usepackage[mathscr]{eucal}
\usepackage{amsmath}
\usepackage{bm} 
\usepackage{amsthm}
\usepackage{amssymb}
\usepackage{amsfonts} 
\usepackage{mathrsfs}
\usepackage{amscd}
\usepackage{enumitem} 
\usepackage{xcolor}
\usepackage{layout}
\usepackage{stmaryrd} 
\usepackage{float} 
\usepackage{longtable}
\usepackage{txfonts} 
\usepackage{setspace} 
\usepackage{subcaption} 

\usepackage{verbatim}
\usepackage{nicefrac}
\usepackage{textcomp}

\usepackage[normalem]{ulem} 
\usepackage{tikz-cd} 

\usepackage{chngcntr} 

\usepackage{hyperref} 

\hypersetup{
    colorlinks=true,
    linkcolor=blue,
    filecolor=blue,      
    urlcolor=blue,
    citecolor = blue,
    }
 
\numberwithin{equation}{section}
\counterwithin{figure}{section}

\newtheorem{theorem}{Theorem}[section]

\newtheorem{definition}[theorem]{Definition}

\theoremstyle{plain}

\newcommand{\vungoc}{V\~u Ng\d{o}c}

\def\epsilon{\varepsilon}
\def\phi{\varphi}
\def\theta{\vartheta}

\newcommand{\ga}{\gamma}

\newcommand{\om}{\omega}

\newcommand{\De}{\Delta}

\def\C{{\mathbb C}}

\def\R{{\mathbb R}}
\def\mbS{{\mathbb S}} 
\def\T{{\mathbb T}}

\def\Z{{\mathbb Z}}

\DeclarePairedDelimiter{\abs}{\lvert}{\rvert}

\DeclareMathOperator{\rk}{rk}

\allowdisplaybreaks

\begin{document}

\title[Hypersemitoric systems: examples and classifications]{Recent examples of hypersemitoric systems and first steps towards a classification: a brief survey}

\author{Yannick Gullentops $\&$ Sonja Hohloch}

\date{\today}

\begin{abstract}
Hypersemitoric systems are 2-degree-of-freedom integrable systems on 4-dimensional symplectic manifolds that have an underlying $\mathbb S^1$-symmetry and no degenerate singularities apart from maybe a finite number of families of so-called parabolic singularities.
We give a short overview of recent examples displaying various bifurcations and sketch a topological-combinatorial classification of the connected components of fibers of hypersemitoric systems.
\end{abstract}


\maketitle



\section{Introduction}

Dynamical systems are a mathematical way to describe changes and processes that appear throughout nature and sciences.
Hamiltonian systems are a more specific class of dynamical systems which appears naturally in many different contexts.
In particular, Hamiltonian systems always have at least one conserved quantity, namely their `energy', but there is a good number of systems that have additional conserved quantities and/or symmetries.
A Hamiltonian systems is said to be {\em completely integrable} if it has the maximal number of independent conserved quantities.

Many classical Hamiltonian systems are completely integrable, such as the spherical pendulum, coupled angular momenta, Gelfand-Zeitlin systems, C. Neumann system, and the Lagrange, Euler, and Kovalevskaya spinning tops.
Moreover, most of these examples have an underlying Hamiltonian $\mbS^1$-action, showing that this is a natural and not uncommon feature.

The aim of this short note is to give a brief motivation for the study of so-called hypersemitoric systems (see Definition \ref{def:hypersemitoric}) and then to summarize the recent results of Gullentops $\&$ Hohloch \cite{gullentopsHohlochCreating,gullentopsHohlochClassifying}.

\subsection*{Acknowledgements}
The first author was fully supported by the UA BOF DocPro4 grant with UA Antigoon ProjectID 34812 and the FWO-FNRS Excellence of Science project G0H4518N with UA Antigoon ProjectID 36584 and the second author partially by both grants.
Both authors wish to thank Joseph Palmer for useful comments.


\section{Toric, semitoric, hypersemitoric systems, and Hamiltonian $\mbS^1$-spaces}

Since this short note focuses on systems on 4-dimensional manifolds, we formulate everything specialized to dimension four:
Given a $4$-dimensional symplectic manifold $(M,\om)$ and a smooth function $g: M \to \R$, the {\em Hamiltonian vector field} $X^g$ is defined via $\om( X^g, \ ) = -dg$ and $g$ is referred to as its {\em Hamiltonian}. Hamilton's equations are given by  the system of (autonomous) ordinary differential equations $z'=X^g(z)$ which can be written in canonical coordinates $(q,p)=(q_1, q_2, p_1, p_2 ) \in \R^2 \times \R^2 \simeq \R^{4}$ with symplectic form $\om = \sum_{i=1}^2 dp_i \wedge dq_i $ as $q'_i= \partial_{p_i}g(q,p)$ and $p'_i= - \partial_{q_i}g(q,p)$ for $1 \leq i \leq 2$.
Using the convention $\mbS^1=\R/2\pi \Z$, having a periodic flow of {\em minimal period $2\pi$} is the same as having an {\em induced effective $\mbS^1$-action}.


\begin{definition}
\label{def:HamSspace}
$(M,\om,J)$ is a {\em Hamiltonian $\mbS^1$-space} if $(M,\om)$ is a compact connected $4$-dimensional symplectic manifold with smooth Hamiltonian $J:M \to \R$ of which the flow is periodic with minimal period $2\pi$.
\end{definition}

Hamiltonian $\mbS^1$-spaces were introduced by Karshon \cite{karshon} and classified up to equivariant symplectomorphism in terms of a labeled graph encoding information about the fixed points and isotropy groups (so-called $\Z_k$-spheres) of the $\mbS^1$-action.


\begin{definition}
\label{def:integrable}
$(M,\om,F=(F_1, F_2))$ is a {\em $4$-dimensional completely (Liouville) integrable system}, briefly an \emph{integrable system}, if $(M,\om)$ is a $4$-dimensional symplectic manifold and $F=(F_1, F_2)\colon M\to\R^2$, referred to as {\em momentum map}, is a smooth map such that $\{F_1, F_2\}:=\om (X^{F_1},X^{F_2})=0$ and $X^{F_1}(p)$ and $ X^{F_2}(p)$ are linearly independent (equivalently, $\rk dF(p)=2)$) for almost all $p\in M$.
\end{definition}

Denote by $\Phi^{F_1}_{t_1}$ and $\Phi^{F_2}_{t_2}$ the Hamiltonian flows of $F_1$ and $F_2$. If they are defined for all times then $F=(F_1, F_2)$ induces an $\R^2$-action via
$$\R^2 \to M, \quad (t, p):=((t_1, t_2), p) \mapsto \Phi^F_t(p):= \bigl( \Phi^{F_1}_{t_1} \circ \Phi^{F_2}_{t_2} \bigr) (p).$$
Points on the manifold are said to be {\em regular} for $F$ if $\rk dF=2$ and otherwise {\em singular}. The value $r \in F(M)$ is {\em regular} if all points in the {\em fiber} $F^{-1}(r)$ are regular. Otherwise $r$ and $F^{-1}(r)$ are said to be {\em singular}. A {\em leaf} is a connected component of a fiber and the {\em leaf space} of $(M, \om, F)$ is obtained by collapsing each leaf to a point. If all fibers are connected, the image $F(M)$ can be identified with the leaf space.

There are symplectic local normal forms (see for example Bolsinov $\&$ Fomenko \cite{bolsinovFomenko} and the references therein): first, near compact connected regular fibers, the so-called {\em action-angle coordinates} provided by the Arnold-Liouville theorem allow to see the system as a nice, trivial $2$-torus bundle with linear flow. Second, near so-called {\em nondegenerate} singular points the system of $F=(F_1, F_2)$ splits into two components which either come from the list {\em regular, elliptic, hyperbolic} or both form a {\em focus-focus} pair.






\begin{definition}
\label{def:toric}
A $4$-dimensional integrable system $(M,\om,F=(F_1, F_2))$ is \emph{toric} if $F$ induces an effective $\T^2 = \mbS^1 \times \mbS^1$-action.
\end{definition}

According to Delzant \cite{delzant}, each compact toric system $(M, \om, F)$ is classified up to equivariant symplectomorphism by the momentum map image $F(M)$ which is a convex polytope with some special properties.
Karshon $\&$ Lerman \cite{karshonLerman} generalized this to the non-compact case.
Note that toric systems only admit singular points with elliptic and/or regular components.

Karshon \cite{karshon} found the necessary and sufficient conditions under which a Hamiltonian $\mbS^1$-space $(M,\om,J)$ extends to a toric integrable system $(M,\om,F=(F_1, F_2))$ with $F_1=J$.

Given a {\em fixed} compact symplectic four manifold, the $\mbS^1$- resp.\ $2$-torus actions on that manifold have been classified by Holm $\&$ Kessler \cite{holmKessler} resp.\ Karshon $\&$ Kessler $\&$ Pinsonnault \cite{karshonKesslerPinsonnault}.


\begin{definition}
\label{def:semitoric}
A $4$-dimensional integrable system $(M,\om,F=(F_1, F_2))$ is \emph{semitoric} if $F_1$ is proper and generates an $\mbS^1$-action of minimal period $2 \pi$ and if all singular points of $F$ are non-degenerate and do not include hyperbolic components. Traditionally, the components of a semitoric system are denoted by $F_1=:J$ and $F_2 =: H$, thus $F=(J, H)$.
\end{definition}

Semitoric systems were introduced by \vungoc\ \cite{vungocPolytope}. They can be seen as generalization of toric integrable systems by dropping the periodicity requirement on the second component. Apart from the singularities appearing in toric systems (regular, elliptic-regular, elliptic-elliptic), semitoric systems may have in addition focus-focus points.
A global symplectic classification was achieved by Pelayo $\&$ \vungoc\ \cite{pelayoVungoc2009, pelayoVungoc2011} together with Palmer $\&$ Pelayo $\&$ Tang \cite{palmerPelayoTang}.
The semitoric classification is valid on compact as well as non-compact manifolds. For the coupled spin oscillators and the coupled angular momenta semitoric systems, all semitoric invariants are by now computed, see Alonso $\&$ Dullin $\&$ Hohloch \cite{adh19, adh20} and the references therein.
Both systems only have maximally one focus-focus point.
The forthcoming work by Alonso $\&$ Hohloch $\&$ Palmer \cite{alonsoHohlochPalmer} together with the article by Alonso $\&$ Hohloch \cite{alonsoHohlochHeight} compute the semitoric invariants of a certain family of semitoric systems with two focus-focus points (which turned out to be significantly harder than for the systems with only one).

Essential examples were constructed in Hohloch $\&$ Palmer \cite{hohlochPalmer2FF}, De Meulenaere $\&$ Hohloch \cite{demeulenaereHohloch}, and Le Floch $\&$ Palmer \cite{leflochPalmer1}. Quantum aspects have been studied in Le Floch $\&$ \vungoc\ \cite{leflochVuNgoc}. Links towards b-symplectic geometry have been established in Brugu\`es $\&$ Hohloch $\&$ Mir $\&$ Miranda \cite{bruguesHohlochMirMiranda}.

Semitoric systems have a natural underlying Hamiltonian $\mbS^1$-action by forgetting $H$.
The relationship between the semitoric classification and Karshon's classification of Hamiltonian $\mbS^1$-spaces is explained by Hohloch $\&$ Sabatini $\&$ Sepe \cite{hss} and an ongoing project by Hohloch $\&$ Sepe $\&$ Sabatini $\&$ Symington.

Recent overviews of the research around semitoric systems and further developments are given for example by Alonso $\&$ Hohloch \cite{alonsoHohlochSurvey}, Pelayo \cite{pelayo}, and Henriksen $\&$ Hohloch $\&$ Martynchuk \cite{henriksenHohlochMartynchuk}.


\begin{definition}
\label{def:hypersemitoric}
 A $4$-dimensional integrable system $(M,\om,F=(F_1, F_2))$ is \emph{hypersemitoric} if $F_1$ is proper and generates an $\mbS^1$-action of minimal period $2 \pi$ and if all degenerate singular points of $F$ (if any) are of parabolic type. Traditionally, the components of a hypersemitoric system are denoted by $F_1=:J$ and $F_2 =: H$, thus $F=(J, H)$.
\end{definition}

Hypersemitoric systems have been first defined in Hohloch $\&$ Palmer \cite{hohlochPalmerLifting}, generalizing semitoric systems by admitting hyperbolic components in singularities. But since the existence of a global $\mathbb S^1$-action prevents the occurrence of hyperbolic-hyperbolic singularities, actually only hyperbolic-regular and/or hyperbolic-elliptic ones may appear in addition to the singularities already occurring in semitoric systems. Since hyperbolic-regular points often give rise to so-called {\em parabolic} (also called {\em cuspidal}) degenerate points it is natural to admit them in the definition of hypersemitoric systems. Parabolic/cuspidal points are named for the `local cusp shape' of their fibers. They can also be seen as `somewhat least degenerate among degenerate points' and are `generic' in the sense that they persist well under perturbations.
For more details, we refer to Efstathiou $\&$ Sugny \cite{efstathiouSugny},  Efstathiou $\&$ Giacobbe \cite{efstathiouGiacobbe}, Bolsinov $\&$ Guglielmi $ \&$ Kudryavtseva \cite{bolsinovGuglielmiKudr}, and Kudryavtseva $\&$ Martynchuk \cite{kudryavtsevaMartynchukCircleAction, kudryavtsevaMartynchukSymplInvar}.

The admittance of hyperbolic components implies that the fibers of hypersemitoric systems are not necessarily connected as this was the case for semitoric systems where the leaf space could be identified with the image of the momentum map.

Hypersemitoric systems have a natural underlying Hamiltonian $\mbS^1$-action induced by $J$. Hohloch $\&$ Palmer \cite{hohlochPalmerLifting} showed how to extend any given Hamiltonian $\mbS^1$-space $(M, \om, J)$ to a hypersemitoric system $(M, \om, F=(J,H))$. Le Floch $\&$ Palmer \cite{leflochPalmer2} study bifurcation theory and other aspects of (hyper)semitoric systems.
An overview of the motivation and research around (hyper)semitoric systems can be found in Henriksen $\&$ Hohloch $\&$ Martynchuk \cite{henriksenHohlochMartynchuk}.


Completing the line of generalisations from toric via semitoric to hypersemitoric systems, we arrive at

\begin{definition}
A 4-dimensional completely integrable system $(M,\omega, F = (F_1, F_2))$ is a {\em proper $\mathbb{S}^1$-system} if $F_1: M \to \R$ is proper and has a periodic flow of minimal period $2 \pi$. Traditionally, the components of a proper $\mathbb{S}^1$-system are denoted by $F_1=:J$ and $F_2 =: H$, thus $F=(J, H)$.
\end{definition}


\section{Explicit examples of hypersemitoric systems}

This subsection summarizes the main results from Gullentops $\&$ Hohloch \cite{gullentopsHohlochCreating} who perturb the explicit (semi)toric system from De Meulenare $\&$ Hohloch \cite{demeulenaereHohloch} to obtain and study a family of hypersemitoric systems which display interesting bifurcations like flaps and swallowtails. The occuring hyperbolic-regular fibres are visualised and show the existence of so-called $k$-stacked tori for $k \in \{2, 3,4\}$.

Now let us be more precise: Consider the standard octagon $\De \subset \R^2$ spanned by the vertices $(0,2)$, $ (0,1)$, $ (1, 0)$, $ (2, 0)$, $ (3, 1)$, $ (3, 2)$, $(2, 3)$ and $ (1, 3)$ which is in fact a so-called {\em Delzant polytope}, i.e., it allows to apply Delzant's \cite{delzant} construction method to obtain a toric system $(M, \om, F=(J,H))$ with underlying $4$-dimensional connected compact symplectic manifold $(M_, \om)$ and $F(M)=\De$. The construction of $M$ starts with symplectic reduction from $\C^8$ so that the points in $M$ can be written as equivalence classes of the form $[z]=[z_1, \dots, z_8]$ where $z_k = x_k + i y_k \in \C$ for $1 \leq k \leq 8$.
The momentum map is then given by
$$
F=(J, H): M \to \R^2 \quad \mbox{with} \quad J([z_1, \dots, z_8])= \frac{1}{2}\abs{z_1}^2, \quad H([z_1, \dots, z_8])= \frac{1}{2} \abs{z_3}^2
$$
and the system $(M, \om, F=(J, H))$ is toric.
For details, we refer the reader to De Meulenare $\&$ Hohloch \cite{demeulenaereHohloch}.
The idea in Gullentops $\&$ Hohloch \cite{gullentopsHohlochCreating} is to perturb the second component of the system $F=(J, H)$ in a suitable way while keeping the first one unchanged. Herefore take parameters $t:=(t_1, t_2, t_3, t_4) \in \R^4$, denote by $\Re$ the real part of a complex number, and consider the functions $\ga_1, \dots, \ga_4: M \to \R$ given by
\begin{align*}
    \gamma_1([z])  := \frac{1}{100} \ \bigl(\overline{z_2z_3z_4}z_6z_7z_8+z_2z_3z_4\overline{z_6z_7z_8}\bigr) = \frac{1}{50} \ \Re(\overline{z_2z_3z_4}z_6z_7z_8)
\end{align*}
and
 \begin{align*}
     \gamma_{2}([z])  : =\frac{1}{50} \ |z_5|^4|z_4|^4, \quad
    \gamma_{3}([z])  : = \frac{1}{50} \ |z_4|^4|z_7|^4,\quad
    \gamma_{4}([z])    : = \frac{1}{100} \ |z_5|^4|z_7|^4.
\end{align*}
Then set
\begin{align*}
   H_t: M \to \R, \quad  H_t:= H_{(t_1,t_2,t_3,t_4)} := (1-2 t_1)H + \sum_{j = 1}^4 t_j\gamma_j \quad \mbox{and} \quad F_t:=(J, H_t): M \to \R^2.
\end{align*}
The perturbation term $\gamma_1$ was (up to a scaling factor) already used by De Meulenaere $\&$ Hohloch \cite{demeulenaereHohloch}. Using in fact four perturbation terms $\ga_1$, $\ga_2$, $\ga_3$, $\ga_4$ is inspired by Le Floch $\&$ Palmer \cite{leflochPalmer1}.

\begin{theorem}[{Gullentops $\&$ Hohloch \cite{gullentopsHohlochCreating}}]
The system $(M, \om, F_t=(J, H_t))$ is completely integrable for all $t \in \R^4$.
The four points in $M$ given by
\begin{align*}
& \left[\sqrt{2}, \ 0, \ 0,\sqrt{2}, \ 2, \ 2\sqrt{2},\sqrt{6},\sqrt{6}\right], && \left[2, \ 2\sqrt{2},\sqrt{6},\sqrt{6},\sqrt{2}, \ 0, \ 0,\sqrt{2} \right], \\
 & \left[\sqrt{2},\sqrt{6},\sqrt{6}, \ 2\sqrt{2}, \ 2,\sqrt{2}, \ 0, \ 0 \right],  && \left[2,\sqrt{2}, \ 0, \ 0,\sqrt{2},\sqrt{6},\sqrt{6}, \ 2\sqrt{2}\right]
\end{align*} 
are singular of rank zero for $F_t=(J, H_t)$ for all $t \in \R^4$. These points lie in $J^{-1}(1) \cup\ J^{-1}(2)$. Any other rank zero points of $F_t=(J, H_t)$ can only appear in $J^{-1}(0) \cup\ J^{-1}(3)$.
Rank one singular points are determined by solving a certain polynomial equation.
\end{theorem}

Note that $J$ still has a periodic flow of minimal period $2 \pi$ since it comes from the toric system $F=(J,H)$. The system $(M, \om, F_t=(J, H_t))$ is for many parameter values $t$ in fact hypersemitoric.

In this context, transitions between elliptic-elliptic and focus-focus singularities turned out to be {\em Hamiltonian Hopf bifurcations}. Roughly, such a bifurcation is {\em supercritical}, if the focus-focus point emerges from an elliptic-elliptic point (think of it for instance as `the elliptic-elliptic value moves from the boundary of the bifurcation diagram into its interior'), and {\em subcritical} if the focus-focus point turns into an `elliptic-elliptic point with a flap'. A rigorous description can be found in Henriksen $\&$ Hohloch $\&$ Martynchuk \cite{henriksenHohlochMartynchuk}.

When hyperbolic-regular singularities occur in the family $F_t=(J, H_t)$, a hyperbolic-regular fiber may consist of a so-called {\em $k$-stacked torus} which can be seen, intuitively, as $k$ tori `placed on top of each other forming a tower of $k$ tori'. In particular, a $2$-stacked torus can be seen as product of a figure-eight loop and a circle.

\begin{theorem}[{Gullentops $\&$ Hohloch \cite{gullentopsHohlochCreating}}]
Given the family of systems $(M, \om, F_t=(J, H_t))$ constructed above, there are explicitly plotted examples for certain parameter values $t$ of hyperbolic singular fibres that are $k$-stacked tori for $k \in \{2, 3,4\}$. Moreover, $k$ can be maximally 13 for the family of systems $(M, \om, F_t=(J, H_t))$.
Twisted tori do not appear. Moreover, there are explicit, visualised examples for flaps and swallowtails and their collisions.
\end{theorem}


\section{Classification of hypersemitoric systems?}

A central point of mathematics is classifications. Concerning symplectic classifications of systems with underlying $\mbS^1$-action, there has been significant progress for the last 35 years: In 1988, Delzant's \cite{delzant} symplectic classification of toric integrable systems; in 1999, Karshon's \cite{karshon} symplectic classification of Hamiltonian $\mbS^1$-spaces; in 2009-2011 and 2019 Pelayo $\&$ \vungoc's \cite{pelayoVungoc2009, pelayoVungoc2011} and Palmer $\&$ Pelayo $\&$ Tang's \cite{palmerPelayoTang} classification of semitoric systems. Symplectic invariants of parabolic orbits and cusp singularities and flaps have been studied by Bolsinov $\&$ Guglielmi $\&$ Kudryavtseva \cite{bolsinovGuglielmiKudr} in 2018 and Kudryavtseva $\&$ Martynchuk \cite{kudryavtsevaMartynchukSymplInvar} in 2021.
Thus the aim clearly is to reach a (symplectic) classification of hypersemitoric systems in the not too distant future.


Under topological aspects, there are already different (semi)local and global classifications of integrable systems like the classification via atoms and molecules by Fomenko, Bolsinov, Oshemkov and many others (see for instance \cite{bolsinovFomenko, bolsinovOshemkov} and the references therein) and the topological one by Zung \cite{zung03}.


In Gullentops \cite{gullentops}, a topological classification of hyperbolic-regular leafs (= connected components of fibers that contain only regular and hyperbolic-regular points) was established. In essence, the idea is to divide the hyperbolic-regular leaf by the $\mbS^1$-action induced by $J$ which leads to a graph with crossings and marked points. Similar thoughts had been independently considered in a different way by Colin de Verdi\`ere $\&$ \vungoc\ \cite{colinDeVerdiereVuNgoc}.

In Gullentops $\&$ Hohloch \cite{gullentopsHohlochClassifying}, these ideas are being generalized to a classification of all types of leaves of a hypersemitoric system by means of a so-called `generalised bouquet' which can be seen as a graph with crossings and marked points which are both in addition decorated with labels coming from the period of the flow of $J$ before passing to the quotient.

This approach, on the one hand, expands and, on the other hand, partially recovers results from Bolsinov $\&$ Fomenko \cite{bolsinovFomenko}, Bolsinov $\&$ Oshemkov \cite{bolsinovOshemkov} and others around the semi-local classification of non-degenerate points of completely integrable systems based on the notion of atoms and from Colin de Verdi\`ere $\&$ \vungoc's \cite{colinDeVerdiereVuNgoc} classification of hyperbolic-regular leafs.

In the symplectic classification of toric systems, the image of the momentum map is the only invariant. For semitoric systems, the image of the momentum map is used to construct the so-called {\em semitoric polygon invariant}. For hypersemitoric systems or even more general $\mbS^1$-systems, there is not yet any symplectic classification in the spirit of the ones for toric and semitoric systems.

One of the first tasks will be to understand how to generalize the semitoric polygon invariant to a meaningful ingredient in a classification of hypersemitoric systems. The main problem to overcome is that fibres of hypersemitoric systems (or even more general $\mbS^1$-systems) are not necessarily connected so that the image of the momentum map cannot be identified with the leaf space.

A possible solution is to work instead with the {\em unfolded bifurcation diagram} (see Efstathiou $\&$ Sugny \cite{efstathiouSugny} and Efstathiou $\&$ Giacobbe \cite{efstathiouGiacobbe}) of a hypersemitoric system $(M, \om, F=(J,H))$ which is the path connected topological space $\mathcal{U}$ together with a projection $\tau :\mathcal{U} \rightarrow F(M)$ such that, for all $r \in \R^2$, the number of points in $\tau^{-1}(r)$ equals the number of connected components in $F^{-1}(r)$. Note that if there is a continuous embedding $F(M) \hookrightarrow M$ that is a right inverse of $F$ then the leaf space forms a {\em canonical} unfolded bifurcation diagram. From a visual point of view, unfolded bifurcation diagrams are well suited to display flaps and swallowtails.

How the other semitoric invariants can be generalized and how all of the invariants can be combined for a classification is a more complicated question, hopefully answered in future projects.



\vspace{3mm}

\noindent
Yannick Gullentops \\
Department of Mathematics \\
University of Antwerp \\
Middelheimlaan 1 \\
B-2020 Antwerp, Belgium \\
{\tt yannick.gullentops@uantwerpen.be}

\vspace{3mm}

\noindent
Sonja Hohloch \\
Department of Mathematics \\
University of Antwerp \\
Middelheimlaan 1 \\
B-2020 Antwerp, Belgium \\
{\tt sonja.hohloch@uantwerpen.be} 

\end{document}